\documentclass[11pt,a4]{article}

\def\funk#1{\noindent {\bf #1{}:}}
\def\be{\begin{equation}}
\def\ee{\end{equation}}
\newcommand{\ff}[1]{{\mbox{\boldmath $#1$}}}
\def\a{\alpha}
\def\b{\beta}
\def\e{\epsilon}

\def\={\approx}
\def\x{\ff{x}}
\def\kk#1#2{\frac{#1}{#2}}
\def\s{\quad}
\def\.{\cdot}

\title{Test Problems in Optimization}

\author{Xin-She Yang \\ 
Department of Engineering, University of Cambridge, \\ Cambridge CB2 1PZ, UK}

\date{}

\begin{document}

\maketitle

\begin{abstract}
{\large 
Test functions are important to validate new optimization algorithms and to compare
the performance of various algorithms.  There are many test functions in the
literature, but there is no standard list or set of test functions one has to follow. 
New optimization algorithms should be tested using at least a subset of functions
with diverse properties so as to make sure whether or not the tested algorithm can
solve certain type of optimization efficiently. Here we provide a selected
list of test problems for unconstrained optimization. \\
}

\noindent {\sf Citation detail: }\\ 
X.-S. Yang, Test problems in optimization, in:
{\it Engineering Optimization: An Introduction with Metaheuristic Applications } 
(Eds Xin-She Yang), John Wiley \& Sons, (2010). 

\end{abstract}

\newpage

In order to validate any new optimization algorithm, we have to validate it
against standard test functions so as to compare its performance with
well-established or existing algorithms. There are many test functions, so there is no
standard list or set of test functions one has to follow. However, various test functions
do exist, so new algorithms should be tested using at least a subset of functions
with diverse properties so as to make sure whether or not the tested algorithm can
solve certain type of optimization efficiently.

In this appendix, we will provide a subset of commonly used test functions with
simple bounds as constraints, though they are often listed as unconstrained problems in literature.
We will list the function form $f(\x)$, its search domain, optimal solutions $\x_*$
and/or optimal objective value $f_*$. Here, we use $\x=(x_1, ..., x_n)^T$
where $n$ is the dimension. \\

%%  \section{Unconstrained Problems}

\funk{Ackley's function} \be f(\x) = -20 \exp\Big[-\kk{1}{5} \sqrt{\kk{1}{n}
\sum_{i=1}^n x_i^2} \; \Big] -\exp\Big[\kk{1}{n} \sum_{i=1}^n \cos ( 2 \pi x_i)  \Big] + 20 +e, \ee
where $n=1,2,...$, and $-32.768 \le x_i \le 32.768$ for $i=1,2,...,n$. This function
has the global minimum $f_*=0$ at $\x_*=(0,0,...,0)$. \\

\funk{De Jong's functions}
The simplest of De Jong's functions is the so-called sphere function
\be f(\x) =\sum_{i=1}^n x_i^2, \s
-5.12 \le x_i \le 5.12, \ee
whose global minimum is obviously $f_*=0$ at $(0,0,...,0)$. This function is unimodal and convex.
A related function is the so-called weighted sphere function or hyper-ellipsoid function
\be f(\x) = \sum_{i=1}^n i x_i^2, \s -5.12 \le x_i \le 5.12, \ee
which is also convex and unimodal with a global minimum $f_*=0$ at $\x_*=(0,0,...,0)$.
Another related test function is the sum of different power function
\be f(\x) = \sum_{i=1}^n |x_i|^{i+1}, \s -1 \le x_i \le 1, \ee
which has a global minimum $f_*=0$ at $(0,0,...,0)$.  \\

\funk{Easom's function} \be f(\x)=-\cos(x) \cos(y) \exp \Big[-(x-\pi)^2 + (y-\pi)^2  \Big], \ee
whose global minimum is $f_*=-1$ at $\x_*=(\pi, \pi)$ within $-100 \le x,y \le 100$.
It has many local minima.
Xin-She Yang extended in 2008 this function to $n$ dimensions, and we have
\be f(\x)=-(-1)^n \Big( \prod_{i=1}^n \cos^2(x_i) \Big) \exp \Big[ - \sum_{i=1}^n (x_i-\pi)^2 \Big], \ee
whose global minimum $f_*=-1$ occurs at $\x_*=(\pi, \pi,..., \pi)$. Here
the domain is $-2 \pi \le x_i \le 2 \pi$ where $i=1,2,...,n$.  \\

\funk{Equality-Constrained Function} \be f(\x) =-(\sqrt{n})^n \prod_{i=1}^n x_i, \ee
subject to an equality constraint (a hyper-sphere) \be \sum_{i=1}^n x_i^2=1. \ee
The global minimum $f_*=-1$ of $f(\x)$ occurs at $\x_*(1/\sqrt{n},...,1/\sqrt{n})$
within the domain $0 \le x_i \le 1$ for $i=1,2,...,n$. \\

\funk{Griewank's function} \be f(\x) = \kk{1}{4000} \sum_{i=1}^n x_i^2 - \prod_{i=1}^n \cos (\kk{x_i}{\sqrt{i}}) +1,
\s -600 \le x_i \le 600, \ee
whose global minimum is $f_*=0$ at $\x_*=(0,0,...,0)$. This function is highly multimodal. \\

\funk{Michaelwicz's function} \be f(\x) = - \sum_{i=1}^n \sin (x_i) \. \Big[ \sin ( \kk{i x_i^2}{\pi}) \Big]^{2m}, \ee
where $m=10$, and $0 \le x_i \le \pi$ for $i=1,2,...,n$.
In 2D case, we have
\be f(x,y) = - \sin (x) \sin^{20} (\kk{x^2}{\pi}) - \sin (y) \sin^{20} (\kk{2 y^2}{\pi}), \ee
where $(x,y) \in [0,5] \times [0,5]$. This function has a global minimum
$f_* \= -1.8013$ at $\x_*=(x_*,y_*)=(2.20319, 1.57049)$. \\

\funk{Perm Functions} \be f(\x) = \sum_{j=1}^n \Big\{ \sum_{i=1}^n (i^j + \b) \Big[(\kk{x_i}{i})^j-1\Big] \Big\}, \s (\b>0), \ee
which has the
global minimum $f_*=0$ at $\x_*=(1,2,...,n)$ in the search domain $-n \le x_i \le n$ for $i=1, ..., n$.
A related function
\be f(\x)=\sum_{j=1}^n \Big\{ \sum_{i=1}^n (i+\b) \Big[x_i^j -(\kk{1}{i})^j \Big] \Big\}^2, \ee
has the global minimum $f_*=0$ at $(1, 1/2, 1/3, ..., 1/n)$ within the bounds $-1 \le x_i \le 1$ for
all $i=1,2,...,n$. As $\b>0$ becomes smaller, the global minimum becomes almost indistinguishable
from their local minima. In fact, in the extreme case $\b=0$, every solution is also a global minimum. \\

\funk{Rastrigin's function} \be f(\x) = 10 n + \sum_{i=1}^n \Big[ x_i^2 - 10 \cos (2 \pi x_i) \Big],
\s -5.12 \le x_i \le 5.12, \ee
whose global minimum is $f_*=0$ at $(0,0,...,0)$. This function is highly multimodal. \\

\funk{Rosenbrock's function} \be f(\x) = \sum_{i=1}^{n-1} \Big[ (x_i-1)^2 + 100 (x_{i+1}-x_i^2)^2 \Big], \ee
whose global minimum $f_*=0$ occurs at $\x_*=(1,1,...,1)$ in the domain
$-5 \le x_i \le 5$ where $i=1,2,...,n$. In the 2D case, it is often written as
\be f(x,y)=(x-1)^2 + 100 (y-x^2)^2, \ee
which is often referred to as the banana function.  \\

\funk{Schwefel's function} \be f(\x) = - \sum_{i=1}^n x_i \sin \Big(\sqrt{|x_i|} \Big),
\s -500 \le x_i \le 500, \ee
whose global minimum $f_*\=-418.9829 n$ occurs at $x_i=420.9687$ where $i=1,2,...,n$. \\

\funk{Six-hump camel back function} \be f(x,y)=(4-2.1x^2 +\kk{1}{3} x^4) x^2 + xy + 4 ( y^2-1) y^2, \ee
where $-3 \le x \le 3$ and $-2 \le y \le 2$. This function has two global minima $f_* \= -1.0316$
at $(x_*,y_*)=(0.0898,-0.7126)$ and $(-0.0898,0.7126)$. \\

\funk{Shubert's function} \be f(\x)=\Big[ \sum_{i=1}^n i \cos \Big( i + (i+1) x \Big) \Big] \.
\Big[ \sum_{i=1}^n i \cos \Big(i + (i+1) y \Big) \Big], \ee
which has 18 global minima $f_* \= -186.7309$ for $n=5$ in the search domain $-10 \le x,y \le 10$. \\

\funk{Xin-She Yang's functions} \be f(\x) =\Big( \sum_{i=1}^n |x_i| \Big) \exp \Big[ -\sum_{i=1}^n \sin (x_i^2) \Big], \ee
which has the global minimum $f_*=0$ at $\x_*=(0,0,...,0)$ in the domain $-2 \pi \le x_i \le 2 \pi$ where $i=1,2,...,n$.
This function is not smooth,  and its derivatives are not well defined at the optimum $(0,0,...,0)$.

A related  function is \be f(\x)=-\Big( \sum_{i=1}^n |x_i| \Big) \exp \Big( -\sum_{i=1}^n x_i^2 \Big),
\s -10 \le x_i \le 10, \ee
which has multiple global minima. For example, for $n=2$, we have $4$ equal minima
$f_* =-1/\sqrt{e} \= -0.6065$ at $(1/2,1/2)$, $(1/2,-1/2)$, $(-1/2,1/2)$ and $(-1/2,-1/2)$.

Yang also designed a standing-wave function with a defect
\be f(\x)=\Big[ e^{-\sum_{i=1}^n (x_i/\b)^{2m}} - 2 e^{-\sum_{i=1}^n x_i^2} \Big]
\. \prod_{i=1}^n \cos^2 x_i, \s m=5, \ee
which has many local minima and the unique global minimum $f_*=-1$ at $\x_*=(0,0,...,0)$
for $\b=15$ within the domain $-20 \le x_i \le 20$ for $i=1,2,...,n$.

The location of the defect can easily be shift to other positions. For example, 
\be f(\x)=\Big[ e^{-\sum_{i=1}^n (x_i/\b)^{2m}} - 2 e^{-\sum_{i=1}^n (x_i-\pi)^2} \Big]
\. \prod_{i=1}^n \cos^2 x_i, \s m=5, \ee
has a unique global minimum $f_*=-1$ at $\x_*=(\pi,\pi,...,\pi)$

Yang also proposed another multimodal function
\be f(\x) = \Big\{\Big[ \sum_{i=1}^n \sin^2(x_i)\Big] - \exp(-\sum_{i=1}^n x_i^2) \Big\} \. \exp\Big[
-\sum_{i=1}^n \sin^2 \sqrt{|x_i|} \; \Big], \ee
whose global minimum $f_*=-1$ occurs at $\x_*=(0,0,...,0)$ in the domain
$-10 \le x_i \le 10 $ where $i=1,2,...,n$. In the 2D case, its landscape looks
like a wonderful candlestick.

Most test functions are deterministic. Yang designed a test function
with stochastic components
\be f(x,y)=-5 e^{-\b  [(x-\pi)^2+(y-\pi)^2]}
-\sum_{j=1}^K \sum_{i=1}^K \e_{ij} e^{-\a  [(x-i)^2 + (y-j)^2]}, \ee
where $\a,\b>0$ are scaling parameters, which can often be taken as $\a=\b=1$.
Here $\e_{ij}$ are random variables and can be drawn from
a uniform distribution $\e_{ij} \sim $ Unif[0,1].
The domain is $0 \le x,y \le K$ and $K=10$. This function has
$K^2$ local valleys at grid locations and the  fixed global
minimum at $\x_*=(\pi,\pi)$. It is worth pointing that
the minimum $f_{\min}$ is random, rather than a fixed value;
it may vary from $-(K^2+5)$ to $-5$, depending $\a$ and $\b$ as well as the
random numbers drawn.

Furthermore, he also designed a stochastic function
\be f(\x) = \sum_{i=1}^n \e_i \; \Big|x_i - \kk{1}{i}\Big|, \s -5 \le x_i \le 5, \ee
where $\e_i \;(i=1,2,...,n)$ are random variables which are uniformly distributed
in $[0,1]$. That is, $\e_i \sim $Unif$[0,1]$.
This function has the unique minimum $f_*=0$ at $\x_*=(1,1/2,...,1/n)$
which is also singular.  \\

\funk{Zakharov's functions} \be f(\x) = \sum_{i=1}^n x_i^2 + \Big (  \kk{1}{2} \sum_{i=1}^n i x_i \Big)^2
+ \Big( \kk{1}{2} \sum_{i=1}^n i x_i \Big)^4, \ee
whose global minimum $f_*=0$ occurs at $\x_*=(0,0,...,0)$.
Obviously, we can generalize this function as
\be f(\x)= \sum_{i=1}^n x_i^2 + \sum_{k=1}^{K} J_n^{2k}, \ee
where $K=1,2,...,20$ and
\be  J_n =\kk{1}{2} \sum_{i=1}^n i x_i. \ee

%% \section{Constrained Problem}


\begin{thebibliography}{50}

\bibitem{Ackley} D. H. Ackley, {\it A Connectionist Machine for Genetic Hillclimbing},
 Kluwer Academic Publishers, 1987.


\bibitem{Floudas} C. A. Floudas, P. M., Pardalos, C. S. Adjiman, W. R. Esposito,
Z. H. Gumus, S. T. Harding, J. L. Klepeis, C. A., Meyer, C. A. Scheiger, {\it Handbook of Test Problems
in Local and Global Optimization}, Springer, 1999.

\bibitem{Hedar} A. Hedar, Test function web pages, http://www-optima.amp.i.kyoto-u.ac.jp /member/student/hedar/Hedar$\_$files/TestGO$\_$files/Page364.htm

\bibitem{Molga} M. Molga, C. Smutnicki, ``Test functions for optimization needs'', \\
http://www.zsd.ict.pwr.wroc.pl/files/docs/functions.pdf

\bibitem{Yang} X.-S. Yang, ``Firefly algorithm, L\'evy flights and global
optimization'', in: {\it Research and Development in Intelligent Systems XXVI},
(Eds M. Bramer et al.), Springer, London, pp. 209-218 (2010).

\bibitem{Yang2} X.-S. Yang and S. Deb, ``Engineering optimization by cuckoo search'', {\it Int. J.
Math. Modeling and Numerical Optimization}, {\bf 1}, No. 4, 330-343 (2010).

\bibitem{Yang3} X.-S. Yang, ``Firefly algorithm, stochastic test functions and design optimization'',
{\it Int. J. Bio-inspired Computation}, {\bf 2}, No. 2, 78-84 (2010).

\end{thebibliography}
\end{document}